\documentclass[12pt,a4paper]{article}
\usepackage{amsfonts}
\usepackage{amssymb}
\usepackage{mathrsfs}
\usepackage{amsmath}
\usepackage{amsthm}
\usepackage{enumerate}
\usepackage{cite}
\usepackage[all]{xy}
\usepackage{extarrows}
\usepackage{indentfirst}
\allowdisplaybreaks
\newtheorem{defn}{Definition}[section]
\newtheorem{thm}[defn]{Theorem}
\newtheorem{lem}[defn]{Lemma}
\newtheorem{prop}[defn]{Proposition}
\newtheorem{cor}[defn]{Corollary}
\newtheorem{eg}[defn]{Example}
\newtheorem{re}[defn]{Remark}

\newcommand\relphantom[1]{\mathrel{\phantom{#1}}}
\newcommand{\bdefn}{\begin{defn}}
\newcommand{\edefn}{\end{defn}}
\newcommand{\bthm}{\begin{thm}}
\newcommand{\ethm}{\end{thm}}
\newcommand{\blem}{\begin{lem}}
\newcommand{\elem}{\end{lem}}
\newcommand{\bprop}{\begin{prop}}
\newcommand{\eprop}{\end{prop}}
\newcommand{\bcor}{\begin{cor}}
\newcommand{\ecor}{\end{cor}}
\newcommand{\beg}{\begin{eg}}
\newcommand{\eeg}{\end{eg}}
\newcommand{\bre}{\begin{re}}
\newcommand{\ere}{\end{re}}
\newcommand{\bpf}{\begin{proof}}
\newcommand{\epf}{\end{proof}}
\newcommand{\benu}{\begin{enumerate}}
\newcommand{\eenu}{\end{enumerate}}
\newcommand{\bc}{\begin{center}}
\newcommand{\ec}{\end{center}}
\newcommand{\bea}{\begin{eqnarray}}
\newcommand{\eea}{\end{eqnarray}}
\newcommand{\Bea}{\begin{eqnarray*}}
\newcommand{\Eea}{\end{eqnarray*}}
\newcommand{\beq}{\begin{equation}}
\newcommand{\eeq}{\end{equation}}
\newcommand{\Beq}{\begin{equation*}}
\newcommand{\Eeq}{\end{equation*}}
\newcommand{\bspl}{\begin{split}}
\newcommand{\espl}{\end{split}}

\begin{document}

\title{\textbf{On  split  Regular Hom-Leibniz algebras}
\author{ Yan Cao$^{1,2},$  Liangyun Chen$^{1}$
 \date{{\small {$^1$ School of Mathematics and Statistics, Northeast Normal
 University,\\
Changchun 130024, China}\\{\small {$^2$  Department of Basic
 Education,
 Harbin University of
Science and Technology,\\ Rongcheng Campus,  Rongcheng 264300,
China}}}}}} \maketitle
\date{}

\begin{abstract}

We introduce the class of split regular Hom-Leibniz algebras as the natural generalization of  split  Leibniz algebras and split regular Hom-Lie algebras. By developing techniques of connections of roots for this kind of algebras,  we show that  such a split regular Hom-Leibniz algebra  $L$ is of the form
$L = U + \sum\limits_{[j] \in \Lambda/\sim}I_{[j]}$ with $U$ a subspace of the abelian  subalgebra $H$ and any $I_{[j]}$, a well described
ideal of $L$, satisfying $[I_{[j]}, I_{[k]}] = 0$ if $[j]\neq [k]$. Under certain conditions, in the case of $L$ being of maximal length,  the simplicity of the algebra is characterized. \\

\end{abstract}
\renewcommand{\thefootnote}{\fnsymbol{footnote}}
\footnote[0]{ Corresponding author(L. Chen): chenly640@nenu.edu.cn.}
\footnote[0]{Supported by  NNSF of China (Nos. 11171055 and
11471090),  NSF of Jilin province (No. 201115006), Scientific
Research Fund of Heilongjiang Provincial Education Department
 (No. 12541184). }

\section{Introduction}
The notion of Hom-Lie algebras was introduced by Hartwig, Larsson and Silvestrov to describe the $q$-deformation of the Witt and the Virasoro algebras \cite {1}. Since then, many authors have studied Hom-type algebras  \cite{2,3,4,5,6,7}. The notion of Leibniz algebras was introduced by Loday \cite{L2}, which is a ``nonantisymmetric'' generalization of Lie algebras. So far, many results of this kind of algebras have been considered in the frameworks of low dimensional
algebras, nilpotence and related problems \cite{AAO, AAO2, AO, B3, B4}. In particular, Makhlouf and Silvestrov introduced the notion of  Hom-Leibniz algebras in \cite {8}, which is a natural generalization of Leibniz algebras and Hom-Lie algebras.

As is well-known, the class of the split algebras is specially related to addition quantum
numbers, graded contractions, and deformations. For instance, for a physical system which displays
a symmetry of $L$, it is interesting to know in detail the structure of the split decomposition
because its roots can be seen as certain eigenvalues which are the additive quantum numbers
characterizing the state of such system. Determining the structure of split algebras will become more and more meaningful in the area of research in mathematical physics. Recently, in \cite{BL52, BL52567, BL5234, BL528}, the structure of arbitrary split Lie algebras, arbitrary split
  Leibniz algebras, arbitrary split  Lie triple systems and arbitrary split regular Hom-Lie algebras have been determined by the techniques of connections of roots.
 The purpose of this paper is to consider the structure of split regular Hom-Leibniz algebras by the techniques of connections of roots based on some work in \cite{BL5234, BL52567}.

Throughout this paper, split regular Hom-Leibniz algebras $L$ are considered of
arbitrary dimension and over an arbitrary base field $\mathbb{K}$.  This paper
is organized as follows. In section 2, we establish the preliminaries on
split regular Hom-Leibniz algebras theory.  In section 3, we show that such
an arbitrary regular Hom-Leibniz algebra $L$ with a symmetric root system is
of the form  $L=U+\sum_{[j]\in \Lambda/\sim} I_{[j]}$ with $U$ a
subspace of the abelian  subalgebra $H$  and any $I_{[j]}$ a well described ideal of
$L$,  satisfying $[I_{[j]},I_{[k]}]
=0$ if $[j]\neq [k]$. In section 4, we show that under certain conditions, in
the case of $L$ being of maximal length,  the simplicity
of the algebra is characterized.

\section{Preliminaries}
First we recall the definitions of  Leibniz algebras, Hom-Lie algebras and Hom-Leibniz algebras.
\bdefn{\rm\cite{L2}}  A \textbf{right Leibniz algebra} $L$ is a vector space over a base field $\mathbb{K}$ endowed with a bilinear product
$[\cdot,\cdot]$ satisfying the Leibniz identity
$$[[y, z], x] = [[y, x], z] + [y, [z, x]],$$
for all $x, y, z \in L$.
\edefn

\bdefn{\rm\cite{5}}  A \textbf{Hom-Lie algebra} $L$ is a vector space over a base field $\mathbb{K}$ endowed with a bilinear product
$$[\cdot,\cdot]:L \times L \rightarrow L$$
and with a linear map $\phi:L \rightarrow L$ such that

$\rm 1$. $[x,y]=-[y,x]$,

$\rm 2$. $[[x, y], \phi(z)] +[[y, z], \phi(x)] + [[z, x], \phi(y)]=0,$

\noindent for all $x, y, z \in L$. When $\phi$ furthermore is an algebra automorphism it is said that $L$ is a \textbf{regular Hom-Lie algebra}.
\edefn

\bdefn{\rm\cite{8}}  A \textbf{Hom-Leibniz algebra} $L$ is a vector space over a base field $\mathbb{K}$ endowed with a bilinear product
$$[\cdot,\cdot]:L \times L \rightarrow L$$
and with a linear map $\phi:L \rightarrow L$ satisfying the Hom-Leibniz identity
$$[[y, z], \phi(x)] = [[y, x], \phi(z)] + [\phi(y), [z, x]],$$
for all $x, y, z \in L$. When $\phi$ furthermore is an algebra automorphism it is said that $L$ is a \textbf{regular Hom-Leibniz algebra}.
\edefn

Clearly Hom-Lie algebras and Leibniz algebras are examples of Hom-Leibniz algebras.

Throughout this paper we will consider regular Hom-Leibniz algebras $L$ being of arbitrary dimension and arbitrary base field $\mathbb{K}$.  $\mathbb{N}$ denotes  the set of all non-negative integers and  $\mathbb{Z}$ denotes the set of all integers.

For any $x\in L$, we consider the adjoint mapping $\mathrm{ad}_{x}:L\rightarrow L$ defined by $\mathrm{ad}_{x}(z)=[z,x].$
A subalgebra $A$ of $L$ is a linear subspace such that $[A,A] \subset A$ and $\phi(A) =
A$. A linear subspace $I$ of $L$ is called an ideal if $[I,L]+[L,I] \subset I$ and $\phi(I) = I$.

Let $L$ be a Hom-Leibniz algebra, the ideal $J(L)$ generated by $\{[x, x]: x \in L\}$ plays an important role in the
theory since it determines the  non-Lie character of $L$. For convenience, write $J$ for   $J(L)$. From the Hom-Leibniz identity, this ideal
satisfies
\begin{equation}\label{188}
[L, J]=0.
\end{equation}

Let us introduce the class of split algebras in the framework of regular Hom-Leibniz
algebras. Denote by $H$ a maximal abelian subalgebra of a Hom-Leibniz algebra $L$. For a linear
functional
$$\alpha:H\rightarrow \mathbb{K},$$
we define the root space of $L$ $($with respect to $H$$)$ associated to $\alpha$ as the subspace
$$L_{\alpha} = \{v_{\alpha}\in L: [v_{\alpha}, h] = \alpha(h)\phi(v_{\alpha}) \ for \ any \ h \in H\}.$$
The elements $\alpha:H\rightarrow \mathbb{K}$ satisfying $L_{\alpha}\neq 0$ are called roots of $L$ with respect to
$H$. We denote $\Lambda:= \{\alpha \in H^{\ast}
\setminus \{0\}: L_{\alpha}\neq 0\}$.

\bdefn
We say that $L$ is a \textbf{split regular Hom-Leibniz algebra}, with respect
to $H$, if
$$L=H\oplus (\oplus_{\alpha \in \Lambda }L_{\alpha}).$$
We also say that $\Lambda$ is the roots system of $L$.
\edefn

Note that when $\phi = \mathrm{Id}$, the split Leibniz algebras become examples of split
regular Hom-Leibniz algebras. Hence, the present paper extends the results in \cite{BL5234}.
For convenience, the mappings $\phi| _{H}$, $\phi|_{H}^{-1}:H \rightarrow H$ will be denoted by $\phi$ and $\phi^{-1}$ respectively.

 \blem \label{3556000}
For any $\alpha$, $\beta \in \Lambda \cup \{0\}$, the following assertions hold.

$\rm 1$. $\phi(L_{\alpha})\subset L_{\alpha \phi^{-1}}$ and $\phi^{-1}(L_{\alpha})\subset L_{\alpha \phi}$.

$\rm 2$. $[L_{\alpha}, L_{\beta}] \subset L_{\alpha\phi^{-1}+\beta\phi^{-1}}$.
\elem

\bpf $\rm 1$.  For $h \in H$ write $h^{'} = \phi(h)$. Then for all $h \in H$ and $v_{\alpha} \in L_{\alpha}$, since
$[v_{\alpha}, h] = \alpha(h)\phi(v_{\alpha})$, one has
$$[\phi(v_{\alpha}),h^{'}] = \phi([v_{\alpha}, h]) = \alpha(h)\phi(\phi(v_{\alpha})) =\alpha\phi^{-1}(h^{'})\phi(\phi(v_{\alpha})).$$
Therefore we get $\phi(v_{\alpha}) \in L_{\alpha\phi^{-1}}$ and so $\phi(L_{\alpha}) \subset L_{\alpha \phi^{-1}}$. In a similar way, one gets
$\phi^{-1}(L_{\alpha})\subset L_{\alpha \phi}$.

$\rm 2$. For any $h \in H$, $v_{\alpha} \in L_{\alpha}$ and $v_{\beta} \in L_{\beta}$, by denoting $h^{'} = \phi(h)$, by Hom-Leibniz identity, we
have that \begin{align*}
[[v_{\alpha}, v_{\beta}], h^{'}]&=[[v_{\alpha}, v_{\beta}], \phi(h)]= [[v_{\alpha},h],\phi(v_{\beta})]+[\phi(v_{\alpha}),[v_{\beta},h]]\\
&=[\alpha(h)\phi(v_{\alpha}),\phi(v_{\beta})]+[\phi(v_{\alpha}), \beta(h)\phi(v_{\beta})]\\
&= (\alpha+\beta)(h)\phi([v_{\alpha}, v_{\beta}])\\
&= (\alpha+\beta)\phi^{-1}(h^{'})\phi([v_{\alpha}, v_{\beta}]).
\end{align*}
Therefore we get $[v_{\alpha}, v_{\beta}] \in L_{\alpha\phi^{-1}+\beta\phi^{-1}}$ and  so  $[L_{\alpha}, L_{\beta}] \subset L_{\alpha\phi^{-1}+\beta\phi^{-1}}$.
\epf

 \blem \label{3556000asd}
The following assertions hold.

$\rm 1$. If $\alpha \in \Lambda$ then $\alpha\phi^{-z} \in \Lambda$ for any $z\in \mathbb{Z}.$

$\rm 2$. $L_{0} = H$.
\elem

\bpf $\rm 1$. It is a consequence of Lemma \ref{3556000}-1.

$\rm 2$. It is clear that  the root space associated to the zero root satisfies $H \subset L_{0}$. Conversely, given any
$v_{0} \in L_{0}$ we can write $$v_{0} =h\oplus(\oplus_{i=1}^{n}v_{\alpha_{i}}),$$
 where $h \in H$ and $v_{\alpha_{i}} \in L_{\alpha_{i}}$
for $i = 1,\cdots, n,$
with $\alpha_{i} \neq \alpha_{j}$ if $i \neq j$. Hence $$0 = [h \oplus(\oplus_{i=1}^{n}v_{\alpha_{i}}), h^{'}]=\oplus_{i=1}^{n}\alpha_{i}(h^{'}) \phi(v_{\alpha_{i}}),$$ for any $h^{'}\in H$. Hence Lemma \ref{3556000}-1 and the fact $\alpha_{i}\neq 0$ give us  $v_{\alpha_{i}}=0$ for $i = 1,\cdots, n$. So $v_{0}=h \in H$.
\epf

\bdefn
 A root system $\Lambda$  of a split Hom-Leibniz algebra is called \textbf{symmetric} if it satisfies that
$\alpha \in \Lambda$ implies $-\alpha \in \Lambda$.
\edefn

\section{Decompositions}
In the following, $L$ denotes a split regular Hom-Leibniz algebra with a symmetric root system $\Lambda$ and $L=
H \oplus(\oplus _{\alpha \in \Lambda}L_{\alpha})$ the corresponding root decomposition. Given a linear functional  $\alpha:H\rightarrow \mathbb{K}$, we denote by $-\alpha:H\rightarrow \mathbb{K}$ the element in $H^{\ast}$ defined by $(-\alpha)(h):=-\alpha(h)$ for all $h\in H$. We begin by developing the techniques of  connections of roots
 in this section.

\bdefn\label{egzmmm}
Let $\alpha$ and $\beta$ be two nonzero roots. We shall say that $\alpha$ is
\textbf{connected} to $\beta$  if there exists $\alpha_{1},\cdots,\alpha_{k}\in \Lambda$  such that

If $k=1$, then

$\rm 1$. $\alpha_{1}\in \{a\phi^{-n}:  n\in \mathbb{N}\}\cap \{\pm \beta \phi^{-m}:  m\in \mathbb{N}\}.$

If $k\geq 2$, then

$\rm 1$.  $\alpha_{1}\in \{a\phi^{-n}:  n\in \mathbb{N}\}$.

$\rm 2$. $\alpha_{1}\phi^{-1}+\alpha_{2}\phi^{-1} \in \Lambda$,

\  \quad  $\alpha_{1}\phi^{-2}+\alpha_{2}\phi^{-2}+\alpha_{3}\phi^{-1} \in \Lambda$,

\  \quad  $\alpha_{1}\phi^{-3}+\alpha_{2}\phi^{-3}+\alpha_{3}\phi^{-2}+\alpha_{4}\phi^{-1} \in \Lambda$,

\  \quad $\cdots \cdots \cdots$

\  \quad $\alpha_{1}\phi^{-i}+\alpha_{2}\phi^{-i}+\alpha_{3}\phi^{-i+1}+\cdots+\alpha_{i+1}\phi^{-1}\in \Lambda$,

\  \quad $\cdots \cdots \cdots$

\  \quad $\alpha_{1}\phi^{-k+2}+\alpha_{2}\phi^{-k+2}+\alpha_{3}\phi^{-k+3}+\cdots+\alpha_{i}\phi^{-k+i}+\cdots+\alpha_{k-1}\phi^{-1}\in \Lambda$.

$\rm 3$. $\alpha_{1}\phi^{-k+1}+\alpha_{2}\phi^{-k+1}+\alpha_{3}\phi^{-k+2}+\cdots+\alpha_{i}\phi^{-k+i-1}+\cdots+\alpha_{k}\phi^{-1}\in \{\pm \beta \phi^{-m}: m \in \mathbb{N}\}$.

\noindent We shall also say
that $\{\alpha_{1},\cdots,\alpha_{k}\}$ is a connection from $\alpha$ to $\beta$.
\edefn

Observe that the case $k =1$ in Definition \ref{egzmmm}  is equivalent to the fact
$\beta=\epsilon \alpha\phi^{z}$ for some $z \in \mathbb{Z}$ and $\epsilon \in \{\pm 1\}$.

By straightforward computations, we can easily get the following  proposition:
\bprop \label{6787777fg}
Let $L$ be a split Hom-Leibniz algebra relative to abelian subalgebra $H$ with roots system $\Lambda$ and $J$  an ideal of $L$ satisfying $[L,J]=0$. Then the semi-direct product $\widehat{L}=L\rtimes L/J$ with respect to the operation given by
$$[(a, x+J), (b, y+J)]=([a,y]-[b,x], [x,y]+J)$$
is a Hom-Lie algebra which has weight decomposition relative to abelian subalgebra $H/J\cap H \subset L/J\subset \widehat{L}$, all adjoints are diagonalizable, and the roots system is exactly  $\Lambda:$
 $$\widehat{L}_{\lambda}=L_{\lambda}\oplus (L_{\lambda}/L_{\lambda}\cap J).$$
\eprop

From Proposition \ref{6787777fg}, we can see roots system $\Lambda$ in a Hom-Leibniz algebra is the same in
a Hom-Lie algebra.  By \cite[Lemmas 2.2, 2.3, Proposition 2.4]{BL52567}, we get the following  Lemmas \ref{3556jjsd}, \ref{3556jjsd77777} and Proposition \ref{alp93556jjsd77777}.

 \blem \label{3556jjsd}
For any $\alpha \in \Lambda$, we have that $\alpha\phi^{z_{1}}$ is connected to  $\alpha\phi^{z_{2}}$ for every
$z_{1}, z_{2} \in \mathbb{Z}.$ We also have that $\alpha \phi^{z_{1}}$ is connected to $-\alpha \phi^{z_{2}}$ in case  $-\alpha \phi^{z_{2}}\in \Lambda$.

\elem

\blem \label{3556jjsd77777}
Let $\{\alpha_{1},\cdots,\alpha_{k}\}$ be a connection from   $\alpha$ to  $\beta$. Then the following
assertions hold.

$\rm 1$. Suppose $\alpha_{1}=\alpha \phi^{-n}, n\in \mathbb{N}$. Then for any $r\in \mathbb{N}$  such that $r\geq n$, there
exists a connection $\{\overline{\alpha}_{1},\cdots,\overline{\alpha}_{k}\}$ from $\alpha$ to $\beta$ such that $\overline{\alpha}_{1}=\alpha \phi^{-r}$.

$\rm 2$. Suppose that $\alpha_{1}=\epsilon\beta\phi^{-m}$ in case $k =1$ or
$$\alpha_{1}\phi^{-k+1}+\alpha_{2}\phi^{-k+1}+\alpha_{3}\phi^{-k+2}+\cdots+\alpha_{k}\phi^{-1}= \epsilon\beta \phi^{-m}$$
in case $k \geq 2$, with $m \in \mathbb{N}$ and $\epsilon \in \{\pm 1\}$. Then for any $r\in \mathbb{N}$ such
that $r\geq m$, there exists a connection  $\{\overline{\alpha}_{1},\cdots,\overline{\alpha}_{k}\}$ from $\alpha$ to $\beta$ such that
 $\overline{\alpha}_{1}=\epsilon\beta \phi^{-r}$ in case $k = 1$ or
$$\overline{\alpha}_{1}\phi^{-k+1}+\overline{\alpha}_{2}\phi^{-k+1}+\overline{\alpha}_{3}\phi^{-k+2}+\cdots+\overline{\alpha}_{k}\phi^{-1}= \epsilon\beta \phi^{-r}$$
in case $k \geq 2$.
\elem

\bprop\label{alp93556jjsd77777}
The relation $\sim$ in $\Lambda$, defined by $\alpha \sim \beta$ if and only if $\alpha$ is connected to $\beta$, is of equivalence.

\eprop

 Proposition \ref{alp93556jjsd77777} tells us the connection relation  $\sim$ in $\Lambda$ is an equivalence relation. So we denote by
$$\Lambda/ \sim: =\{[\alpha]:\alpha \in \Lambda\},$$
where $[\alpha]$ denotes the set of nonzero roots of $L$ which are connected to $\alpha$.
Our next goal is to associate an adequate ideal $I_{[\alpha]}$ to any $[\alpha]$. For a fixed $\alpha \in \Lambda$, we define
$$I_{0,[\alpha]}:=\mathrm{span_{\mathbb{K}}}\{[L_{\beta}, L_{-\beta}]: \beta \in [\alpha]\}\subset H$$
and
$$V_{[\alpha]}:= \oplus_{\beta \in [\alpha]}L_{\beta}.$$
Then we denote by $I_{[\alpha]}$ the direct sum of the two subspaces above, that is,
$$I_{[\alpha]}:= I_{0,[\alpha]}\oplus V_{[\alpha]}.$$

\bprop\label{a8107217771}
For any $\alpha \in \Lambda$, the linear subspace $I_{[\alpha]}$ is a subalgebra of $L$.
\eprop
\bpf
First, it is sufficient  to check that $I_{[\alpha]}$ satisfies $[I_{[\alpha]}, I_{[\alpha]}] \subset I_{[\alpha]}$. By $I_{0,[\alpha]}\subset H$,
 it is clear that $[I_{0,[\alpha]}, I_{0,[\alpha]}]$=0 and we  have
 \begin{equation}\label{810227}
[ I_{0,[\alpha]}\oplus V_{[\alpha]},  I_{0,[\beta]}\oplus V_{[\beta]}] \subset [I_{0,[\alpha]}, V_{[\beta]}] + [V_{[\alpha]}, [I_{0,[\beta]}]+[V_{[\alpha]},V_{[\beta]}].
\end{equation}

Let us consider the first summand in (\ref{810227}). For $\beta \in [\alpha]$,  by Lemmas  \ref{3556000} and \ref{3556jjsd}, one gets
$[I_{0,[\alpha]}, L_{\beta}] \subset L_{\beta\phi^{-1}}$, where $\beta\phi^{-1}\in [\alpha]$. Hence
 \begin{equation}\label{810227810227}
[I_{0,[\alpha]}, V_{[\alpha]}]\subset V_{[\alpha]}.
\end{equation}
 Similarly, we can also get
  \begin{equation}\label{810227810227as}
 [V_{[\alpha]}, [I_{0,[\alpha]}]\subset V_{[\alpha]}.
 \end{equation}
Next, we consider  the third summand in   (\ref{810227}). Given
$\beta$,
 $\gamma \in [\alpha]$ such that $[L_{\beta},L_{\gamma}]\neq 0$, if
 $\gamma= -\beta$, we have $[L_{\beta},L_{\gamma}
]=[L_{\beta},L_{-\beta}] \subset I_{0,[\alpha]}.$ Suppose
$\gamma \neq -\beta$, by Lemma \ref{3556000}-2,
one gets $\beta\phi^{-1}+\gamma\phi^{-1}\in \Lambda$. Therfore, we get  $\{\beta, \gamma\}$ is a connection from $\beta$ to
 $\beta\phi^{-1}+\gamma\phi^{-1}$. The transitivity of $\sim$  gives  that $\beta\phi^{-1}+\gamma\phi^{-1}\in [\alpha]$ and so
 $[L_{\beta},L_{\gamma}
]\subset L_{\beta\phi^{-1}+\gamma\phi^{-1}}\subset V_{[\alpha]}$. Hence
$$[\oplus_{\beta \in [\alpha]}L_{\beta},\oplus_{\beta \in [\alpha]}L_{\beta}]\subset I_{0,[\alpha]}\oplus V_{[\alpha]}.$$ That is,
\begin{equation}\label{810721810227}
[V_{[\alpha]},V_{[\alpha]}]\subset I_{[\alpha]}.
\end{equation}
From (\ref{810227}), (\ref{810227810227}), (\ref{810227810227as}) and (\ref{810721810227}), we get $[I_{[\alpha]}, I_{[\alpha]}]\subset I_{[\alpha]}.$

Second, we have to verify that $\phi(I_{[\alpha]})=I_{[\alpha]}$. It is a direct consequence of Lemmas \ref{3556000}-1 and  \ref{3556jjsd}.
\epf

\bprop\label{a810721777}
If $[\alpha] \neq [\beta]$, then $[I_{[\alpha]}, I_{[\beta]}] = 0$.
\eprop
\bpf
We have
 \begin{equation}\label{810227080504}
[ I_{0,[\alpha]}\oplus V_{[\alpha]},  I_{0,[\beta]}\oplus V_{[\beta]}] \subset [I_{0,[\alpha]}, V_{[\beta]}] + [V_{[\alpha]}, [I_{0,[\beta]}]+[V_{[\alpha]},V_{[\beta]}].
\end{equation}
Let us consider the third summand $[V_{[\alpha]},V_{[\beta]}]$ in  (\ref{810227080504}) and suppose there exist $\alpha_{1} \in [\alpha]$ and  $\alpha_{2} \in [\beta]$ such that $[L_{\alpha_{1}}, L_{\alpha_{2}}] \neq 0.$  By known condition $[\alpha] \neq [\beta]$,  one gets $\alpha_{1}\neq -\alpha_{2}$. So $\alpha_{1}\phi^{-1}+\alpha_{2}\phi^{-1}\in \Lambda$.
Hence $\{\alpha_{1},\alpha_{2},-\alpha_{1}\phi^{-1}\}$ is a connection from  $\alpha_{1}$ to  $\alpha_{2}$. By the transitivity of
the connection relation, we have $\alpha \in [\beta]$, a contradiction. Hence $[L_{\alpha_{1}}, L_{\alpha_{2}}]=0$ and so
\begin{equation}\label{2005}
[V_{[\alpha]},V_{[\beta]}]=0.
\end{equation}
Next we consider  the first summand $[I_{0,[\alpha]}, V_{[\beta]}]$ in  (\ref{810227080504}). Let us take $\alpha_{1} \in [\alpha]$ and $\alpha_{2} \in [\beta]$
and conclude that
\begin{equation}\label{bag}
 \alpha_{2} ([L_{\alpha_{1}}, L_{-\alpha_{1}}])=0.
\end{equation}
Indeed, by applying  Hom-Leibniz identity and  (\ref{2005}), one gets
 \begin{equation}\label{ruler}
[\phi(L_{\alpha_{2}}),[L_{\alpha_{1}}, L_{-\alpha_{1}}]] = 0.
\end{equation}
By $\phi$ is an algebra automorphism and (\ref{ruler}), one gets
 \begin{equation}\label{ruler12}
\phi[(L_{\alpha_{2}}),\phi^{-1}[L_{\alpha_{1}}, L_{-\alpha_{1}}]] = 0,
\end{equation}
 that is
 \begin{equation}\label{cup}
[(L_{\alpha_{2}}),\phi^{-1}[L_{\alpha_{1}}, L_{-\alpha_{1}}]] = 0,
\end{equation}
 where $\phi^{-1}[L_{\alpha_{1}}, L_{-\alpha_{1}}]\subset H$.
Hence (\ref{cup}) gives
\begin{equation}\label{cupbag}
\alpha_{2}\phi^{-1}([L_{\alpha_{1}}, L_{-\alpha_{1}}])=0,
\end{equation}
 for any $\alpha_{1} \in [\alpha]$ and $\alpha_{2} \in [\beta]$. By  Lemma \ref{3556000}-1 and $\phi$ is an algebra automorphism, we get
$$\phi([L_{\alpha_{1}}, L_{-\alpha_{1}}])\subset [L_{\alpha_{1}\phi^{-1}}, L_{-\alpha_{1}\phi^{-1}}],$$
that is  $$[L_{\alpha_{1}}, L_{-\alpha_{1}}]\subset \phi^{-1}([L_{\alpha_{1}\phi^{-1}}, L_{-\alpha_{1}\phi^{-1}}])$$ and by (\ref{cupbag}), one gets $$\alpha_{2} ([L_{\alpha_{1}}, L_{-\alpha_{1}}])=0.$$ From  $[L_{\alpha_{2}},[L_{\alpha_{1}}, L_{-\alpha_{1}}]]\subset \alpha_{2} ([L_{\alpha_{1}}, L_{-\alpha_{1}}])\phi(L_{\alpha_{2}})=0$, we prove that $[I_{0,[\alpha]}, V_{[\beta]}]=0$. In a similar way, we get $[V_{[\alpha]}, I_{0,[\beta]}]=0$ and we conclude, together with  (\ref{810227080504}) and (\ref{2005}), that $[I_{[\alpha]}, I_{[\beta]}] = 0$.
 \epf

\bdefn\label{2hi}
A Hom-Leibniz algebra $L$ is said to be \textbf{simple} if its product is nonzero and its only ideals are
$\{0\}$, $J$ and $L$.
\edefn

It should be noted that the above definition agrees with the definition of simple Hom-Lie algebra, since
$J=\{0\}$ in this case.

\bthm\label{cao}
The following assertions hold.

$\rm 1$. For any $\alpha \in \Lambda$, the subalgebra
 $$I_{[\alpha]}=I_{0,[\alpha]}\oplus V_{[\alpha]}$$
of $L$ associated to $[\alpha]$ is an ideal of $L$.

$\rm 2$. If $L$ is simple, then there exists a connection from $\alpha$ to $\beta$ for any $\alpha,\beta \in \Lambda$
and $H=\sum_{\alpha \in \Lambda}[L_{\alpha}, L_{-\alpha}]$.
\ethm

\bpf
$\rm 1$. Since $[I_{[\alpha]}, H]+[H, I_{[\alpha]}]=[I_{[\alpha]}, L_{0}]+[L_{0}, I_{[\alpha]}]\subset V_{[\alpha]}$,
taking into account Propositions \ref{a8107217771} and  \ref{a810721777}, we have
$$[I_{[\alpha]}, L]=[I_{[\alpha]}, H\oplus(\oplus_{\beta \in [\alpha]}L_{\beta})\oplus(\oplus_{\gamma \not \in [\alpha]}L_{\gamma})]\subset I_{[\alpha]}$$
and
$$[L, I_{[\alpha]}]=[H\oplus(\oplus_{\beta \in [\alpha]}L_{\beta})\oplus(\oplus_{\gamma \not \in [\alpha]}L_{\gamma}), I_{[\alpha]}]\subset I_{[\alpha]}.$$
As we also have by  Lemmas \ref{3556000}-1 and  \ref{3556jjsd} that $\phi( I_{[\alpha]})= I_{[\alpha]}$, we conclude that $I_{[\alpha]}$ is an ideal of $I$.

$\rm 2$. The simplicity of $L$ implies $I_{[\alpha]}\in \{J,L\}$ for any $\alpha \in \Lambda$. If some $\alpha \in \Lambda$ is such that $I_{[\alpha]}=L$ then $[\alpha]=\Lambda$. Hence, $L$ has all of its nonzero roots connected and $H=\sum_{\alpha \in \Lambda}[L_{\alpha}, L_{-\alpha}]$.  Otherwise, if $I_{[\alpha]}=J$ for any $\alpha \in \Lambda$ then $[\alpha]=[\beta]$ for any $\alpha, \beta \in \Lambda$ and so $[\alpha]=\Lambda$. We also conclude that $L$ has all of its nonzero roots connected and $H=\sum_{\alpha \in \Lambda}[L_{\alpha}, L_{-\alpha}]$.
\epf

\bthm
For a vector space complement $U$ of $span_{\mathbb{K}}\{{[L_{\alpha},L_{-\alpha}]: \alpha \in \Lambda}\}$ in H, we have
$$L = U + \sum\limits_{[\alpha] \in \Lambda/\sim}I_{[\alpha]},$$
where any $I_{[\alpha]}$ is one of the ideas  of $L$ described in Theorem \ref{cao}-1, satisfying
$[I_{[\alpha]},I_{[\beta]}]=0,$ whenever $[\alpha] \neq [\beta].$
\ethm

\bpf
Each $I_{[\alpha]}$ is well defined and, by Theorem \ref{cao}-1, an ideal of $L$. It is
clear that
$$L=H\oplus (\oplus_{\alpha \in \Lambda}L_{\alpha})=U + \sum\limits_{[\alpha] \in \Lambda/\sim}I_{[\alpha]}.$$
Finally Proposition \ref{a810721777} gives us  $[I_{[\alpha]}, I_{[\beta]}] = 0$ if $[\alpha] \neq [\beta]$.
\epf

\bdefn\label{2hijian}
The \textbf{annihilator} of a Hom-Leibniz algebra $L$ is the set $\mathrm{Z}(L) = \{x \in L: [x, L] + [L, x] = 0\}$.
\edefn

\bcor\label{667}
If $Z(L) = 0$ and $[L,L] = L$, then $L$ is the direct sum of the ideals given in Theorem \ref{cao},
$$L =\oplus_{[\alpha] \in \Lambda/\sim}I_{[\alpha]}.$$
\ecor
\bpf
From $[L,L]=L$, it is clear that $L=\oplus_{[\alpha] \in \Lambda/\sim}I_{[\alpha]}.$ Finally, the sum is direct because $Z(L) = 0$ and $[I_{[\alpha]}, I_{[\beta]}]=0$ if
 $[\alpha] \neq [\beta]$.
\epf

\section{The simplicity of Split regular Hom-Leibniz algebras of maximal length. }

In this section we focus on the simplicity of split regular Hom-Leibniz algebras by centering our attention in those of maximal length. From now on char($\mathbb{K}$)=0.
\bdefn
We say that a split regular Hom-Leibniz algebra $L$ is of \textbf{maximal length} if $\mathrm{dim}L_{\alpha}$=1 for any $\alpha \in \Lambda$.
\edefn

\blem \label{lemma 4.1}
Let $L$ be a split regular Hom-Leibniz algebra with $\mathrm{Z}(L)=0$ and $I$ an ideal of $L$.  If $I\subset H$ then $I=\{0\}$.
\elem
\bpf
Suppose there exists a nonzero ideal $I$ of $L$ such that $I \subset H$. We get $[I, H]+[H, I] \subset [H, H]=0$.  We also get $[I, \oplus_{\alpha \in \Lambda}L_{\alpha}]+[\oplus_{\alpha \in \Lambda}L_{\alpha}, I]\subset I\subset H$. Then taking into account $H = L_{0}$, we have  $[I, \oplus_{\alpha \in \Lambda}L_{\alpha}]+[\oplus_{\alpha \in \Lambda}L_{\alpha}, I]\subset H \cap (\oplus_{\alpha \in \Lambda}L_{\alpha})=0.$ From here $I\subset  \mathrm{Z}(L)=0$, which is a contradiction.
\epf

\blem \label{lemma 4.2}
For any $\alpha, \beta \in \Lambda$  with $\alpha \neq \beta$ there exists $h_{0} \in  H$ such that
$\alpha(h_{0})\neq 0$ and  $\alpha(h_{0})\neq \beta(h_{0})$.
\elem

\bpf This can be proved completely analogously to \cite[Lemma 4.2]{BL52567}.
\epf

\blem \label{lemma 4.3}
Let $L = H\oplus(\oplus_{\alpha \in \Lambda} L_{\alpha})$ be a split regular Hom-Leibniz algebra. If $I$ is an ideal of $L$ then $I=(I\cap H)\oplus(\oplus_{\alpha \in \Lambda}(I\cap L_{\alpha})).$
\elem

\bpf
Let $x \in I$. We can write $x = h+\sum_{j=1}^{n}v_{\alpha_{j}}$ with $h \in H$, $v_{\alpha_{j}} \in L_{\alpha_{j}}$ and $\alpha_{j}\neq \alpha_{k}$ if $j\neq k$. Let us show that any $v_{\alpha_{j}} \in I$.

If $n = 1$ we write $x = h+v_{\alpha_{1}} \in I$. By taking $h^{'} \in H$ such that $\alpha_{1}(h^{'})\neq 0$  we have $[x, h^{'}] = \alpha_{1}(h^{'})\phi(v_{\alpha_{1}}) \in I$ and so $\phi(v_{\alpha_{1}}) \in  I$. Therefore $\phi^{-1}(\phi(v_{\alpha_{1}}))=v_{\alpha_{1}} \in I.$

Suppose now $n > 1$ and consider $\alpha_{1}$ and  $\alpha_{2}$. By Lemma \ref{lemma 4.2} there exists
$h_{0} \in H$ such that $\alpha_{1}(h_{0})\neq 0 $ and $\alpha_{1}(h_{0})\neq  \alpha_{2}(h_{0})$. Then we have
\begin{equation}\label{270}
I\ni [x, h_{0}]=\alpha_{1}(h_{0})\phi(v_{\alpha_{1}})+\alpha_{2}(h_{0})\phi(v_{\alpha_{2}})+\cdots+\alpha_{n}(h_{0})\phi(v_{\alpha_{n}})
\end{equation}
and

\begin{equation}\label{270588}
I\ni \phi(x)=\phi(h)+\phi(v_{\alpha_{1}})+\phi(v_{\alpha_{2}})+\cdots+\phi(v_{\alpha_{n}}).
\end{equation}

\noindent By multiplying (\ref{270588}) by $\alpha_{2}(h_{0})$ and subtracting (\ref{270}), one gets
\begin{align*}
&\alpha_{2}(h_{0})\phi(h)+(\alpha_{2}(h_{0})-\alpha_{1}(h_{0}))\phi(v_{\alpha_{1}})+(\alpha_{2}(h_{0})-\alpha_{3}(h_{0}))\phi(v_{\alpha_{3}})\\
&+\cdots+(\alpha_{2}(h_{0})-\alpha_{n}(h_{0}))\phi(v_{\alpha_{n}})\in I.
\end{align*}

\noindent  By denoting $ \widetilde{h}:= \alpha_{2}(h_{0})\phi(h) \in H$  and $v_{\alpha_{i}\phi^{-1}}:=(\alpha_{2}(h_{0})-\alpha_{i}(h_{0}))\phi(v_{\alpha_{i}})\in L_{\alpha_{i}\phi^{-1}}$, we can write
\begin{equation}\label{2705886004}
 \widetilde{h}+v_{\alpha_{1}\phi^{-1}}+v_{\alpha_{3}\phi^{-1}}+\cdots+v_{\alpha_{n}\phi^{-1}} \in I.
\end{equation}
Now we can argue as above, with (\ref{2705886004}), to get
$$\widetilde{\widetilde{h}}+v_{\alpha_{1}\phi^{-2}}+v_{\alpha_{4}\phi^{-2}}+\cdots+v_{\alpha_{n}\phi^{-2}} \in I,$$
for $\widetilde{\widetilde{h}} \in H$ and $v_{\alpha_{1}\phi^{-2}} \in L_{\alpha_{i}\phi^{-2}}.$ Following this process, we obtain
$$\overline{h}+v_{\alpha_{1}\phi^{-n+1}} \in I,$$
with $\overline{h} \in H$ and $v_{\alpha_{1}\phi^{-n+1}} \in L_{\alpha_{i}\phi^{-n+1}}.$ As in the above case $n=1$, we conclude that $v_{\alpha_{1}\phi^{-n+1}} \in I$
and consequently $v_{\alpha_{1}}=\phi^{-n+1}(v_{\alpha_{1}\phi^{-n+1}}) \in I$.

In a similar way we can prove that $v_{\alpha_{i}} \in I$ for $i \in \{2,\cdots,n\}$, and the proof is complete.
\epf

Let us return to a split regular Hom-Leibniz algebra of maximal length $L$. From now on $L = H\oplus(\oplus_{\alpha \in \Lambda} L_{\alpha})$ denotes  a split Hom-Leibniz algebra of maximal length. By Lemma \ref{lemma 4.3}, we assert that given any nonzero ideal $I$ of $L$ then

\begin{equation}\label{8889765566}
I=(I\cap H)\oplus(\oplus_{\alpha \in \Lambda^{I}} L_{\alpha}),
\end{equation}
where $\Lambda^{I}:=\{\alpha \in  \Lambda: I\cap  L_{\alpha}\neq 0\}$.

In particular, case $I$=$J$, we get
\begin{equation}\label{888976556689}
J=(J\cap H)\oplus(\oplus_{\alpha \in \Lambda^{J}} L_{\alpha}).
\end{equation}

\noindent From here, we can write
\begin{equation}\label{88897655668977}
\Lambda=\Lambda^{J}\cup \Lambda^{\neg J},
\end{equation}

\noindent where
$$\Lambda^{J}:=\{\alpha \in \Lambda: L_{\alpha}\subset J\}$$
\noindent and
$$\Lambda^{\neg J}:=\{\alpha \in \Lambda: L_{\alpha}\cap J=0\}.$$
\noindent Therefore
\begin{equation}\label{8889765566897756}
L = H\oplus(\oplus_{\alpha \in \Lambda^{\neg J}} L_{\alpha})\oplus(\oplus_{\beta \in \Lambda^{J}} L_{\beta}).
\end{equation}

We note that the fact that $L=[L,L]$, the split decomposition given by   (\ref{8889765566897756}) and (\ref{188})
 show

\begin{equation}\label{88897655668977565678}
H = \sum_{\alpha \in \Lambda^{\neg J}}[L_{\alpha},L_{-\alpha}].
\end{equation}

Now, observe that the concept of connectivity of nonzero roots given in Definition \ref{egzmmm} is not strong
enough to detect if a given $\alpha \in \Lambda$ belongs to $\Lambda^{J}$ or to $\Lambda^{\neg J}$. Consequently we lose the information
respect to whether a given root space $L_{\alpha}$ is contained in $J$ or not, which is fundamental to study the
simplicity of $L$. So, we are going to refine the concept of connections of nonzero roots in the setup of
maximal length split  regular Hom-Leibniz algebras. The symmetry of $\Lambda^{J}$  and  $\Lambda^{\neg J}$ will be understood as usual. That
is, $\Lambda^{\gamma}$,$\gamma \in \{J, \neg J\}$, is called symmetric if $\alpha \in \Lambda^{\gamma}$ implies $-\alpha \in \Lambda^{\gamma}$.

\bdefn\label{666}
Let $\alpha, \beta \in \Lambda^{\gamma}$  with $\gamma  \in \{J, \neg J\}$. We say that $\alpha$ is \textbf{$\neg J$-connected} to $\beta$, denoted by
$\alpha\sim _{\neg J} \beta$, if there exist
$$\alpha_{2}, \cdots, \alpha_{k} \in \Lambda^{\neg J}$$
\noindent such that

\noindent $\rm 1$. $\{\alpha_{1},\alpha_{1}\phi^{-1}+\alpha_{2}\phi^{-1}, \alpha_{1}\phi^{-2}+\alpha_{2}\phi^{-2}+\alpha_{3}\phi^{-1},\cdots, \alpha_{1}\phi^{-k+1}+\alpha_{2}\phi^{-k+1}+\alpha_{3}\phi^{-k+2}+\cdots+\alpha_{i}\phi^{-k+i-1}+\cdots+\alpha_{k}\phi^{-1}\}\subset \Lambda^{\gamma},$

\noindent $\rm 2$. $\alpha_{1}\in \alpha\phi^{-n}$, for $n \in \mathbb{N},$

\noindent $\rm 3$. $\alpha_{1}\phi^{-k+1}+\alpha_{2}\phi^{-k+1}+\alpha_{3}\phi^{-k+2}+\cdots+\alpha_{i}\phi^{-k+i-1}+\cdots+\alpha_{k}\phi^{-1}\in \pm\beta \phi^{-m},$ for $m \in \mathbb{N}$.

 We shall also say
that $\{\alpha_{1},\alpha_{2},\cdots,\alpha_{k}\}$ is a $\neg J$-connection from $\alpha$ to $\beta$.
\edefn

\bprop \label{67890345777777}
 The following assertions  hold.

$\rm 1$. If  $\Lambda^{\neg J}$  is symmetric, then the relation $\sim_{\neg J}$ is an equivalence relation in $\Lambda^{\neg J}$.

 $\rm 2$. If $L= [L, L]$ and $\Lambda^{\neg J}$, $\Lambda^{J}$ are symmetric, then the relation $\sim_{\neg J}$ is an equivalence relation in $\Lambda^{J}$.
\eprop

\bpf
$\rm 1$. Can be proved in a similar way to Proposition \ref{alp93556jjsd77777}.

$\rm 2$.  Let $\beta \in \Lambda^{J}$. Since $\beta \neq 0$, (\ref{88897655668977565678}) gives us that there exists $\alpha \in \Lambda^{\neg J}$
 such that $[\phi(L_{\beta}), [L_{\alpha}, L_{-\alpha}]]\neq 0$. By Hom-Leibniz identity, either $[[L_{\beta}, L_{\alpha}], \phi(L_{-\alpha})] \neq 0$ or $[[L_{\beta}, L_{-\alpha}], \phi(L_{\alpha})] \neq 0$. In the first case, the
$\neg J$-connection $\{\beta, \alpha, -\alpha \phi^{-1}\}$ gives us $\beta\sim_{\neg J} \beta$  while in the second one the $\neg J$-connection $\{\beta, -\alpha, \alpha \phi^{-1}\}$
gives us the same result. Consequently  $\sim_{\neg J} $ is reflexive in $\Lambda^{J}$. The symmetric and transitive character
of $\sim_{\neg J} $ in $\Lambda^{J}$  follows as in Proposition \ref{alp93556jjsd77777}.
 \epf

Let us introduce the notion of root-multiplicativity in the framework of split regular Hom-Leibniz algebras of
maximal length, in a similar way to the ones for split regular Hom-Lie algebras $($see
\cite {BL52567} $)$.

\bdefn\label{8IU}
We say that a split regular Hom-Leibniz  algebra of maximal length $L$ is \textbf{root-multiplicative} if the
below conditions hold.

$\rm 1$. Given $\alpha, \beta \in \Lambda^{\neg J}$ such that $\alpha\phi^{-1}+\beta\phi^{-1} \in \Lambda$ then $[L_{\alpha}, L_{\beta}] \neq 0$.

$\rm 2$. Given  $\alpha \in \Lambda^{\neg J}$ and $\gamma \in \Lambda^{J}$ such that  $\alpha\phi^{-1}+\gamma\phi^{-1} \in \Lambda^{J}$ then $[L_{\alpha}, L_{\gamma}] \neq 0$.
\edefn

Another interesting notion related to split regular Hom-Leibniz algebras of maximal length $L$ is those of Lie-annihilator.
Write $L = H\oplus(\oplus_{\alpha \in \Lambda^{\neg J}} L_{\alpha})\oplus(\oplus_{\beta \in \Lambda^{J}} L_{\beta})$
(see  (\ref{8889765566897756})).

\bdefn
The \textbf{Lie-annihilator} of a split Hom-Leibniz algebra of maximal length $L$ is the set
\begin{align*}
\mathrm{Z}_{\mathrm{Lie}}(L)=&\Big\{x \in L: [x, H\oplus(\oplus_{\alpha \in \Lambda^{\neg J}} L_{\alpha})]+[H\oplus(\oplus_{\alpha \in \Lambda^{\neg J}} L_{\alpha}), x]=0\Big\}.
\end{align*}
\edefn

Clearly we  have $\mathrm{Z}(L)\subset \mathrm{Z}_{\mathrm{Lie}}(L)$.

\bprop \label{67890345777777999}
Suppose $L= [L, L]$ and  $L$ is root-multiplicative. If $\Lambda^{\neg J}$
 has all of its roots $\neg J$-connected,
then any ideal $I$ of $L$ such that $I \not \subseteq H\oplus J$ satisfies $I = L$.
\eprop

\bpf
By (\ref{8889765566}) and (\ref{88897655668977}), we can write
$$I=(I\cap H)\oplus(\oplus_{\alpha_{i} \in \Lambda^{\neg J,I}} L_{\alpha_{i}})\oplus(\oplus_{\beta_{j} \in \Lambda^{ J,I}} L_{\beta_{j}}),$$
\noindent where $\Lambda^{\neg J,I}:=\Lambda^{\neg J}\cap \Lambda^{I}$ and $\Lambda^{ J,I}:=\Lambda^{ J}\cap \Lambda^{I}.$  Since $I \not \subseteq H\oplus J$, one gets $
\Lambda^{\neg J,I}\neq \emptyset$ and so we can take some $\alpha_{0} \in \Lambda^{\neg J,I}$ such that
\begin{equation}\label{888976444589}
L_{\alpha_{0}}\subset I.
\end{equation}
 By  Lemma \ref{3556000}-1, $\phi(L_{\alpha_{0}})\subset L_{\alpha_{0}\phi^{-1}}$. Since $L$ is of maximal length, we have  $ 0\neq \phi(L_{\alpha_{0}})= L_{\alpha_{0}\phi^{-1}}$. (\ref{888976444589}) and $\phi$ is injective give us $\phi(L_{\alpha_{0}})\subset \phi(I)=I$. So,  $L_{\alpha_{0}\phi^{-1}}\subset I$. Similarly we get
 \begin{equation}\label{888976444}
 L_{\alpha_{0}\phi^{-n}}\subset I, \ \mathrm{for} \ n \in \mathbb{N}.
  \end{equation}
  For any $\beta \in \Lambda^{\neg J}$, $\beta \not \in \pm \alpha_{0}\phi^{-n}$, for $n \in \mathbb{N}$, the fact that $\alpha_{0}$ and $\beta$ are $\neg J$-connected gives us a  $\neg J$-connection
$\{\gamma_{1},\cdots,\gamma_{k}\}\subset \Lambda^{\neg J}$ from $\alpha_{0}$ to $\beta$ such that

\noindent $\rm 1$. $\{\gamma_{1}, \gamma_{1}\phi^{-1}+\gamma_{2}\phi^{-1}, \gamma_{1}\phi^{-2}+\gamma_{2}\phi^{-2}+\gamma_{3}\phi^{-1},\cdots, \gamma_{1}\phi^{-k+1}+\gamma_{2}\phi^{-k+1}+\gamma_{3}\phi^{-k+2}+\cdots+\gamma_{k}\phi^{-1}\} \subset \Lambda^{\neg J}$,

\noindent $\rm 2$. $\gamma_{1} \in \alpha_{0}\phi^{-n},$ for $n \in \mathbb{N}$,

\noindent $\rm 3$.  $ \gamma_{1}\phi^{-k+1}+\gamma_{2}\phi^{-k+1}+\gamma_{3}\phi^{-k+2}+\cdots+\gamma_{k}\phi^{-1} \in \pm\beta\phi^{-m},$ for $m \in \mathbb{N}$.

Consider $\gamma_{1}, \gamma_{2}$ and $\gamma_{1}\phi^{-1}+\gamma_{2}\phi^{-1}$. Since $\gamma_{1}, \gamma_{2} \in \Lambda^{\neg J},$ the root-multiplicativity and maximal length
of $L$ show $[L_{\gamma_{1}}, L_{\gamma_{2}}]=L_{\gamma_{1}\phi^{-1}+\gamma_{2}\phi^{-1}},$ and by  (\ref{888976444}), $L_{\gamma_{1}}\subset I$. So we have
$$L_{\gamma_{1}\phi^{-1}+\gamma_{2}\phi^{-1}}\subset I.$$
\noindent We can argue in a similar way from $\gamma_{1}\phi^{-1}+\gamma_{2}\phi^{-1}, \gamma_{3}$, and $\gamma_{1}\phi^{-2}+\gamma_{2}\phi^{-2}+\gamma_{3}\phi^{-1}$ to get $$L_{\gamma_{1}\phi^{-2}+\gamma_{2}\phi^{-2}+\gamma_{3}\phi^{-1}}\subset I.$$
\noindent Following this process with the $\neg J$-connection
$\{\gamma_{1},\cdots,\gamma_{k}\}$,  we obtain that
$$L_{\gamma_{1}\phi^{-k+1}+\gamma_{2}\phi^{-k+1}+\gamma_{3}\phi^{-k+2}+\cdots+\gamma_{k}\phi^{-1}}\subset I.$$
\noindent From here, we get that either
\begin{equation}\label{888976444469}
 L_{\beta\phi^{-m}}\subset I\ \mathrm{or} \ L_{-\beta\phi^{-m}}\subset I,
\end{equation}
for any $\beta \in \Lambda^{\neg J}$, $m \in \mathbb{N}$. Note that $\beta \in \Lambda^{\neg J}$ gives us
\begin{equation}\label{opq}
\beta\phi^{-m} \in \Lambda^{\neg J}, \ \mathrm{for} \ m \in \mathbb{N}.
\end{equation}

\noindent Since $H=\sum_{\beta \in \Lambda^{\neg J}}[L_{\beta},L_{-\beta}]$, by (\ref{888976444469}) and (\ref{opq}), we get
\begin{equation}\label{13ff}
H\subset I.
\end{equation}

Now, given any $\delta \in \Lambda$, the facts $\delta\neq 0$, $H\subset I$ and the maximal length of $L$ show that
\begin{equation}\label{13fffggg}
[L_{\delta}, H]=L_{\delta}\subset I.
\end{equation}

\noindent From (\ref{13ff}) and (\ref{13fffggg}), we conclude $I=L$.
\epf

\bprop \label{67891110345777777999}
Suppose $L = [L,L]$, $\mathrm{Z}(L) = 0$ and $L$ is root-multiplicative. If  $\Lambda^{\neg J}$, $\Lambda^{J}$ are symmetric
 and $\Lambda^{ J}$ has all of its roots $\neg J$-connected, then any nonzero ideal $I$ of $L$ such that $I \subseteq J$ satisfies either
$I = J$ or $J = I \oplus K$ with $K$ an ideal of $L$.
\eprop

\bpf
By (\ref{8889765566}) and (\ref{88897655668977}), we can write
$$I=(I\cap H)\oplus(\oplus_{\alpha_{i} \in \Lambda^{ J,I}} L_{\alpha_{i}}),$$
where $\Lambda^{ J,I}\subset \Lambda^{ J}$. Observe that the fact $\mathrm{Z}(L) = 0$ implies
\begin{equation}\label{888778555}
J\cap H=\{0\}.
\end{equation}
Indeed, we have $[L_{\alpha}, J\cap H]+[J\cap H, L_{\alpha}]\subset [L, J]=0$, for any $\alpha \in \Lambda^{ J}$ and $[H, J\cap H]+[J\cap H, H]=0$. So,  $[L, J\cap H]+[J\cap H, L]=0$. That is, $J\cap H\subset \mathrm{Z}(L) = 0$. Hence we can write
 $$I=\oplus_{\alpha_{i} \in \Lambda^{ J,I}} L_{\alpha_{i}},$$
 with $\Lambda^{ J,I}\neq \emptyset$, and so we can take some $\alpha_{0}\in \Lambda^{ J,I}$ such that  $L_{\alpha_{0}}\subset I$. We can argue with the
  root-multiplicativity
and the maximal length of $L$ as in Proposition \ref{67890345777777999} to conclude that given any $\beta \in \Lambda^{ J}$,
there exists a $\neg J$-connection $\{\gamma_{1},\cdots, \gamma_{k}\}$ from $\alpha_{0}$ to $\beta$ such that
$$[[\cdots[L_{\gamma_{1}}, L_{\gamma_{2}}],\cdots], L_{\gamma_{k}}] \in L_{\pm \beta\phi^{-m}}, \ \  \mathrm{for}  \ \ m \in \mathbb{N}$$
and so
\begin{equation}\label{888976444DHGaf44}
 L_{\epsilon \beta \phi^{-m}}\subset I, \ \ \mathrm{for} \ \ \mathrm{some} \ \epsilon \in \pm 1,\ \  m \in \mathbb{N}.
   \end{equation}
Note that $\beta \in \Lambda^{J}$ indicates $L_{\beta}\subset J$. By $\phi$ is injective, we get $\phi(L_{\beta})\subset \phi(J)=J$.  By  Lemma \ref{3556000}-1, $\phi(L_{\beta})\subset L_{\beta\phi^{-1}}$. Since $L$ is of maximal length, we have  $ 0\neq \phi(L_{\beta})= L_{\beta\phi^{-1}}$. So,  $L_{\beta\phi^{-1}}\subset J$.
Similarly we get
 \begin{equation}\label{888976444DHG}
 L_{\beta\phi^{-m}}\subset J, \ \mathrm{for} \ m \in \mathbb{N}.
  \end{equation}
By (\ref{888976444DHGaf44}) and (\ref{888976444DHG}), we easily get
\begin{equation}\label{3347009}
\epsilon_{\beta}\beta \phi^{-m} \in \Lambda^{ J,I} \ \mathrm{for} \ \mathrm{any} \ \beta \in \Lambda^{J},  \  \ \mathrm{some} \ \epsilon_{\beta}\in \pm 1 \ \ \mathrm{and} \ \ m \in \mathbb{N}.
  \end{equation}

  Suppose  $-\alpha_{0}\in \Lambda^{ J,I}$. Then we also have that $\{-\gamma_{1},\cdots, -\gamma_{k}\}$ from $-\alpha_{0}$ to $\beta$  is a $\neg J$-connection from
 $-\alpha_{0}$ to $\beta$ satisfying
 $$[[\cdots[L_{-\gamma_{1}}, L_{-\gamma_{2}}],\cdots], L_{-\gamma_{k}}] \in L_{-\epsilon_{\beta}\beta\phi^{-m}}\subset I$$
 and so $L_{\beta\phi^{-m}}+L_{-\beta\phi^{-m}}\subset I$. Hence, (\ref{888976556689}) and  (\ref{888778555}) imply that $I=J$.

 Now,  suppose there is not any  $\alpha_{0}\in \Lambda^{ J,I}$ such that  $-\alpha_{0}\in \Lambda^{ J,I}$.  (\ref{3347009}) allows us to write $\Lambda^{ J}=
 \Lambda^{ J,I}\cup (-\Lambda^{ J,I})$. By denoting $K = \oplus_{\alpha_{i} \in \Lambda^{ J,I}} L_{-\alpha_{i}}$, we have
 \begin{equation}\label{basket}
J = I \oplus K.
  \end{equation}
 Let us finally show that $K$ is an ideal of $L$. We have $[L, K]\subset [L, J]=0$ and
 $$[K, L]\subset [K, H]+ [K, \oplus_{\beta \in \Lambda^{\neg J}}L_{\beta}]+[K, \oplus_{\gamma \in \Lambda^{J}}L_{\gamma}]\subset K+[K, \oplus_{\beta \in \Lambda^{\neg J}}L_{\beta}].$$
 \noindent Let us consider the last summand $[K, \oplus_{\beta \in \Lambda^{\neg J}}L_{\beta}]$ and suppose there exist $\alpha_{i} \in \Lambda^{ J,I}$ and $\beta \in \Lambda^{\neg J}$ such that $[L_{-\alpha_{i}}, L_{\beta}]\neq 0.$ Since $L_{-\alpha_{i}}\subset K \subset J$, we get $-\alpha_{i}\phi^{-1}+\beta\phi^{-1} \in \Lambda^{J}.$
 By the root-multiplicativity of $L$,
the symmetries of $\Lambda^{ J}$ and $\Lambda^{\neg J}$, and the fact $L_{\alpha_{i}}\subset I$, one gets  $0 \neq [L_{\alpha_{i}}, L_{-\beta}]=L_{\alpha_{i}\phi^{-1}-\beta\phi^{-1}}\subset I$, that is, $\alpha_{i}\phi^{-1}-\beta\phi^{-1} \in \Lambda^{ J,I}$.  Hence, $-\alpha_{i}\phi^{-1}+\beta\phi^{-1} \in -\Lambda^{ J,I}$
and so
$[L_{-\alpha_{i}}, L_{\beta}]\subset K$. Consequently $[K, \oplus_{\beta \in \Lambda^{\neg J}}L_{\beta}]\subset K$.

Next, we have to verify that $\phi(K)=K$. Indeed, since $I$, $J$ are two nonzero ideals, we have $\phi(I)=I$ and $\phi(J)=J$. By (\ref{basket}) and $\phi$ is an algebra automorphism, it is easily to get
 $\phi(K)=K$. We conclude that $K$ is an ideal of $L$.
\epf

We introduce the definition of primeness in the framework of Hom-Leibniz algebras following the same
motivation that in the case of simplicity (see Definition \ref{2hi}  and the above paragraph).

\bdefn
A Hom-Leibniz algebra $L$ is said to be \textbf{prime} if given two ideals $I$, $K$ of $L$ satisfying $[I, K] +
[K, I]=0$, then either $I \in \{0, J, L\}$ or $K \in \{0, J, L\}$.
\edefn

We also note that the above definition agrees with the definition of prime Hom-Lie algebra, since $J= {0}$
in this case.

Under the hypotheses of Proposition \ref{67891110345777777999} we have:

\bcor\label{555}
If furthermore $L$ is prime, then any nonzero ideal $I$ of $L$ such that $I\subseteq J$ satisfies $I= J$.
\ecor
\bpf
Observe that, by Proposition \ref{67891110345777777999}, we could have $J = I \oplus K$ with $I, K$  ideals of $L$, being $[I, K] +
[K, I] = 0$ as consequence of $I, K \subseteq J$. The primeness of $L$ completes the proof.
\epf

\bprop \label{67877999}
Suppose $L= [L, L]$, $\mathrm{Z}_{\mathrm{Lie}} (L)=0$ and $L$ is root-multiplicative. If $\Lambda^{\neg J}$
 has all of its roots $\neg J$-connected, then any ideal $I$ of $L$ such that $I \not \subseteq  J$ satisfies $I = L$.
\eprop
\bpf
Taking into account Lemma \ref{lemma 4.3} and Proposition \ref{67890345777777999} we just have to study the case in which
$$I=(I\cap H)\oplus(\oplus_{\beta_{j} \in \Lambda^{ J,I}} L_{\beta_{j}}),$$
\noindent where $I\cap H\neq 0$. But this possibility never happens. Indeed, observe that
 $[L_{\alpha}, I\cap H]+[I\cap H, L_{\alpha}]\subset [L_{\alpha},  H]+[ H, L_{\alpha}]\subset  L_{\alpha}$, for any $\alpha \in \Lambda^{\neg J}$ and $[L_{\alpha}, I\cap H]+[I\cap H, L_{\alpha}]\subset [L_{\alpha}, I]+[I, L_{\alpha}]\subset I$. So,  $[L_{\alpha}, I\cap H]+[I\cap H, L_{\alpha}]\subset L_{\alpha}\cap I=0$, for $\alpha \in \Lambda^{\neg J}$. Since we also have $[I\cap H,H]+[H, I\cap H]\subset [ H, H]=0$, one gets $I\cap H \subset  \mathrm{Z}_{\mathrm{Lie}}(L)=0$, a contradiction.  Proposition \ref{67890345777777999} completes the proof.
\epf

Given any $\alpha \in \Lambda^{\gamma}$, $\gamma \in \{J, \neg J\}$, we denote by
$$\Lambda_{\alpha}^{\gamma}:= \{\beta \in  \Lambda^{\gamma}: \beta \sim_{\neg J} \alpha\}.$$
For $\alpha \in  \Lambda^{\gamma}$, we write $H_{\Lambda_{\alpha}^{\gamma}}:=\mathrm{span}_{\mathbb{K}}\{[ L_{\beta},L_{-\beta}]: \beta \in \Lambda_{\alpha}^{\gamma}\}\subset H,$ and $V_{\Lambda_{\alpha}^{\gamma}}:=\oplus_{\beta \in \Lambda_{\alpha}^{\gamma}}L_{\beta}$. We denote by $L_{\Lambda_{\alpha}^{\gamma}}$ the following subspace of $L$, $L_{\Lambda_{\alpha}^{\gamma}}:=H_{\Lambda_{\alpha}^{\gamma}}\oplus V_{\Lambda_{\alpha}^{\gamma}}.$

\blem\label{2222}
If $L=[L,L]$, then $L_{\Lambda_{\alpha}^{J}}$ is an ideal of $L$ for any $\alpha \in \Lambda^{J}$.
\elem

\bpf By (\ref{88897655668977565678}), we get $H_{\Lambda_{\alpha}^{J}}=0$ and so $$L_{\Lambda_{\alpha}^{J}}=\oplus_{\beta \in \Lambda_{\alpha}^{J}}L_{\beta}.$$

\noindent It is easy to see that
\begin{equation}\label{basket456}
[L_{\delta},L_{\Lambda_{\alpha}^{J}}]+ [L_{\Lambda_{\alpha}^{J}},L_{\delta}]\subset [L, J]=0,\ \  \mathrm{for} \ \ \delta \in \Lambda^{J}
 \end{equation}
and
\begin{equation}\label{basket456er}
[L_{\Lambda_{\alpha}^{J}},H]+[H, L_{\Lambda_{\alpha}^{J}}]\subset L_{\Lambda_{\alpha}^{J}}.
  \end{equation}
 We will show that
  \begin{equation}\label{basket456fgj}
  [L_{\Lambda_{\alpha}^{J}},L_{\gamma}]\subset L_{\Lambda_{\alpha}^{J}},\  \ \mathrm{for} \ \  \mathrm{any} \ \ \gamma \in  \Lambda^{\neg J}.
   \end{equation}
  Indeed, given any $\beta \in \Lambda_{\alpha}^{J}$ such that $[L_{\beta}, L_{\gamma}]\neq 0$ we have $\beta\phi^{-1}+\gamma\phi^{-1} \in \Lambda^{J}$ and so $\{\beta, \gamma\}$ is a $\neg J$-connection from $\beta$ to $\beta\phi^{-1}+\gamma\phi^{-1}$. By the symmetry and transitivity of $\sim_{\neg J}$ in $\Lambda^{J}$,   we get $\beta\phi^{-1}+\gamma\phi^{-1} \in  \Lambda_{\alpha}^{J}$. Hence $[L_{\beta},L_{\gamma}]\subset L_{\Lambda_{\alpha}^{J}}$, that is, (\ref{basket456fgj}) holds.  taking into account (\ref{8889765566897756}), by (\ref{basket456}), (\ref{basket456er}) and (\ref{basket456fgj}), we have
    \begin{equation}\label{basket456fgjFULL}
   [L, L_{\Lambda_{\alpha}^{J}}]+[L_{\Lambda_{\alpha}^{J}}, L]\subset L_{\Lambda_{\alpha}^{J}}.
  \end{equation}

 Next, we have to verify that
   \begin{equation}\label{ttt6fgjFULL}
 \phi(L_{\Lambda_{\alpha}^{J}})=L_{\Lambda_{\alpha}^{J}}.
  \end{equation}
  Indeed, given  $\beta \in \Lambda_{\alpha}^{J}$ such that  $L_{\beta}\subset L_{\Lambda_{\alpha}^{J}}$ and we have $[L_{\beta}, J]=0$. By  $\phi$ is an algebra automorphism, one gets
   \begin{equation}\label{basket456fgjFULL345}
   [\phi(L_{\beta}), \phi(J)]= [\phi(L_{\beta}), J]=0.
 \end{equation}
 By Lemma \ref{3556000}-1 and $L$ is of maximal length, we get $ 0 \neq \phi(L_{\beta})=L_{\beta\phi^{-1}}.$ Therefore, (\ref{basket456fgjFULL345}) gives us $L_{\beta\phi^{-1}} \in \Lambda^{J}$. That is,  $\phi(L_{\beta})\subset L_{\tau},\   \mathrm{for }\  \tau \in  \Lambda^{J}$. Since  $L$ is of maximal length,
  $\phi(L_{\beta})= L_{\tau},\   \mathrm{for }\  \tau \in  \Lambda^{J}$. Hence (\ref{ttt6fgjFULL}) holds. Consequently, it follows from (\ref{basket456fgjFULL}) and (\ref{ttt6fgjFULL}) that
 $L_{\Lambda_{\alpha}^{J}}$ is an ideal of $L$.
\epf

\bthm
Suppose $L = [L, L]$, $\mathrm{Z}_{\mathrm{Lie}}(L) = 0$, $L$ is root-multiplicative. If $\Lambda^{J}$, $\Lambda^{\neg J}$ are symmetric
then $L$ is simple if and only if it is prime and $\Lambda^{J}$, $\Lambda^{\neg J}$ have all of their roots $\neg J$-connected.
\ethm

\bpf
Suppose $L$ is simple. If $\Lambda^{J}\neq \emptyset$ and we take $\alpha \in \Lambda^{J}$.  Lemma \ref{2222} gives us $L_{\Lambda_{\alpha}^{J}}$
is a nonzero ideal
of $L$.  By $L$ is simple, one gets $L_{\Lambda_{\alpha}^{J}}
= J = \oplus_{\beta \in \Lambda^{J}} L_{\beta}$ (see  (\ref{888976556689}) and (\ref{888778555})). Hence, $\Lambda_{\alpha}^{J}=\Lambda^{J}$
 and consequently $\Lambda^{J}$ has all of its roots $\neg J$-connected.

Consider now any $\gamma \in \Lambda^{\neg J}$ and the subspace $L_{\Lambda_{\gamma}^{\neg J}}$.
Let us denote by $I(L_{\Lambda_{\gamma}^{\neg J}})$ the ideal of $L$ generated by $L_{\Lambda_{\gamma}^{\neg J}}$.
We observe that the fact $J$ is an ideal of $L$ and we assert that $I(L_{\Lambda_{\gamma}^{\neg J}})\cap (\oplus_{\delta \in \Lambda^{\neg J}}L_{\delta})$ is contained
in the linear span of the set
\begin{align*}
 &\{[[\cdots[v_{\gamma^{'}},v_{\alpha_{1}}],\cdots], v_{\alpha_{n}}];\ \  [v_{\alpha_{n}},[\cdots[v_{\alpha_{1}}, v_{\gamma^{'}}],\cdots]];\\
 &[[\cdots[v_{\alpha_{1}},v_{\gamma^{'}}],\cdots], v_{\alpha_{n}}]; \ \ [v_{\alpha_{n}},[\cdots[v_{\gamma^{'}}, v_{\alpha_{1}}],\cdots]];\\
 & \mathrm{where} \ 0\neq v_{\gamma^{'}} \in L_{\Lambda_{\gamma}^{\neg J}}, \ 0\neq v_{\alpha_{i}}\in L_{\alpha_{i}},\  \alpha_{i} \in \Lambda^{\neg J} \mathrm{and}\  n \in \mathbb{N}\}.
 \end{align*}
By simplicity $I(L_{\Lambda_{\gamma}^{\neg J}})=L$. From here, given any $\delta \in \Lambda^{\neg J}$, the above observation gives us that we can write $\delta=\gamma^{'}\phi^{-m}+\alpha_{1}\phi^{-m}+\alpha_{2}\phi^{-m+1}+\cdots+\alpha_{m}\phi^{-1}$ for $\gamma^{'} \in \Lambda_{\gamma}^{\neg J}$, $\alpha_{i} \in \Lambda^{\neg J}$, $m \in \mathbb{N}$ and being the partial sums nonzero. Hence $\{\gamma^{'},\alpha_{1},\cdots,\alpha_{m}\}$ is a $\neg J$-connection from $\gamma^{'}$ to $\delta$. By the symmetry and transitivity of $\sim_{\neg J}$ in $\Lambda^{\neg J}$, we deduce $\gamma$ is $\neg J$-connected to any $\delta \in  \Lambda^{\neg J}$. Consequently, Proposition \ref{67890345777777} gives us that
$\Lambda^{\neg J}$ has all of its roots $\neg J$-connected.
Finally, since $L$ is simple then is prime.

The converse is a consequence of Corollary \ref{555}  and Proposition \ref
{67877999}.
\epf

\end{document}